\documentclass[runningheads]{iwap}

\usepackage{amsmath,amssymb}
\usepackage{url}
\usepackage{graphicx} 
\usepackage{natbib}
\usepackage[]{hyperref}
\newcommand{\Nat}{\mathbb N}

\newcommand{\F}{{\cal F}}
\newcommand{\T}{{\cal T}}

\newcommand{\M}{\mathcal{M}}

\newcommand{\ti}{\widetilde}
\newcommand{\Ind}{{\mbox{\sc I}}}

\newcommand*{\esssup}{\mathop{\mathrm{ess\,sup}}\displaylimits}
\def\esssup{\operatornamewithlimits{ess\,sup}}

\begin{document}
%
\title*{Time management in a Poisson fishing model
}
%
\toctitle{Time management in a Poisson fishing model
}
%
\titlerunning{A double stopping in a fishing model
}
%
\author{
  Anna Karpowicz\inst{1}\inst{2}
  \and 
 Krzysztof Szajowski\inst{1}\inst{3}
}
%
\index{Szajowski, K.}
\index{Karpowicz, A.}
%
\authorrunning{Karpowicz and Szajowski} 

%
\institute{
  Institute of Mathematics and Computer Science\\
 Wroc\l{}aw University of Technology\\
 Wyb\-rze\-\.{z}e Wys\-pia\'{n}\-skie\-go~27\\
 PL-50-370~Wroc\-\l{}aw, Poland
   \and
  (e-mail: {anna.karpowicz@pwr.wroc.pl})
  \and
  (e-mail: {krzysztof.szajowski@pwr.wroc.pl})
}

\maketitle             

\begin{abstract}
The aim of the paper is to extend the model of "fishing problem". The simple formulation is following. The angler goes to fishing. He buys fishing ticket for a fixed time. There are two places for fishing at the lake. The fishes are caught according to renewal processes which are different at both places. The fishes' weights and the inter-arrival times are given by the sequences of i.i.d. random variables with known distribution functions. These distributions are different for the first and second fishing place. The angler's satisfaction measure is given by difference between the utility function dependent on size of the caught fishes and the cost function connected with time. On each place the angler has another utility functions and another cost functions. In this way, the angler's relative opinion about these two places is modeled. For example, on the one place better sort of fish can be caught with bigger probability or one of the places is more comfortable. Obviously our angler wants to have as much satisfaction as possible and additionally he have to leave the lake before the fixed moment. Therefore his goal is to find two optimal stopping times in order to maximize his satisfaction. The first time corresponds to the moment, when he eventually should change the place and the second time, when he should stop fishing. These stopping times have to be less than the fixed time of fishing.
The value of the problem and the optimal stopping times are derived.

\end{abstract}

\section{Introduction}
The solution of double optimal stopping problem, in so called "fishing problem", will be presented. One of the first author who considered the basic version of this problem was ~\cite{sta74:captime} and further generalizations were done by ~\cite{stawoo74:gone}, ~\cite{stawarwoo76:delobs}, ~\cite{krasta90:sizedep}. The detailed review of the papers connected with "fishing problem" was presented by ~\cite{fer97:pois}. The simple formulation of our double optimal stopping problem is following. The angler goes to fishing. He buys fishing ticket for a fixed time $t_0$. There are two places for fishing at the lake and he can change the place at any moment $s$. The fishes are caught according to renewal processes $\{N_i(t),t\geq 0\}$, where $N_i(t)$ is the number of claims during the time $t$ at the place $i \in \{1,2\}$. Let $T_{i,n}$ denote the moment of $n$-th claim at the place $i$ (we fix $T_{1,0}=0$ and $T_{2,0}=s$), then the random variables $S_{i,n}=T_{i,n}-T_{i,n-1}$ are i.i.d. with continuous distribution function (c.d.f. for short) $F_i$.  The fishes' weights, which were caught at the place $i$, are given by the sequence of i.i.d. random variables $X_{i,0},X_{i,1},X_{i,2},\dots$ with c.d.f. $H_i$ (we fix $X_{1,0}=0$ and $X_{2,0}=0$). The renewal process is independent on the sequence of claims. The angler's satisfaction measure is given by difference between the utility function $g_i: [0,\infty) \rightarrow [0,G_i]$ dependent on size of the caught fishes and the cost function $c_i: [0,t_0] \rightarrow [0,C_i] $ connected with time. We assume that $g_i$ and $c_i$ are continuous
and bounded, additionally $c_i$ is differentiable. On each place the angler has different utility functions and cost functions. In this way, the angler's relative opinion about these two places is modeled. For example, on the one position better sort of fish can be caught with bigger probability or one of the piers is more comfortable. The angler can change the place of fishing at any time $s$. The mass of the fishes $M_t^s$ which were caught up to time $t$ if the change of the position took place at the time $s$ ($M_t=M_t^t$) is given by
	\begin{equation*}
		M_t^s=\sum_{n=1}^{N_1(s\wedge t)}X_{1,n} + \sum_{n=1}^{N_2((t-s)^+)}X_{2,n}.
	\end{equation*}
Let $Z(s,t)$ denote the angler's payoff for stopping at time $t$ if the change of the position took place at time $s$. The payoff can be expressed as 
	\begin{equation}\label{AK Z 1}
		Z(s,t)=\left\{
			\begin{array}{ll}
 				g_1(M_t)-c_1(t)               &\mbox{ if } t < s \leq t_0, \\
 				g_1(M_s)-c_1(s)+g_2(M_t^s-M_s)-c_2(t-s) &\mbox{ if } s \leq t \leq t_0,\\
 				- C &\mbox{ if }  t_0<t,
			\end{array}\right.
	\end{equation} 
where $C=C_1+C_2$. With the notation
$w_2(m,s,\ti m,t)=w_1(m,s)+g_2(\ti m - m)-c_2(t-s)$ and $w_1(m,t)=g_1(m)-c_1(t)$, formula~(\ref{AK Z 1}) reduces to
	\begin{equation}\label{Z(s,t)}
		Z(s,t)=\Ind_{\{t < s \leq t_0\}}w_1(M_t,t)+\Ind_{\{s \leq t \leq t_0\}}w_2(M_s,s,M_t^s,t)-\Ind_{\{t_0<t\}}C
	\end{equation}

\section{The optimization problem}
Let $\F_t=\sigma(X_{1,0},T_{1,0},X_{1,1},T_{1,1},\dots,X_{1,N_1(t)},T_{1,N_1(t)})$ be the $\sigma$-field generated by all events up to time $t$, if there was no change of parameters
and $\F_{s,t}=\sigma(\F_s,X_{2,0},T_{2,0},\dots,X_{2,N_2(t-s)},T_{2,N_2(t-s)})$, the $\sigma$-field generated by all events up to time $t$ if the change of parameters was at time $s$. 
For simplicity of notation we set $\F_n:=\F_{T_{1,n}}$, $\F_{s,n}:=\F_{s,T_{2,n}}$.  
Let $\M(\F_n)$ denote the set of non-negative and $\F_n$-measurable random variables. From now on, $\T$ and $\T^s$ stands for the sets of stopping times with respect to the $\sigma$-fields $\{\F_t, t\geq 0\}$
and $\{\F_{s,t}, 0\leq s \leq t\}$, respectively. Furthermore, define for $\ n \in \Nat$ and $n\leq K$ the sets
		$\T_{n,K}=\{\tau \in \T:\tau \geq 0,\ T_{1,n}\leq \tau \leq T_{1,K}\}$ 
and
		$\T_{n,K}^s=\{\tau \in \T^s: 0 \leq s \leq \tau,\ T_{2,n}\leq \tau \leq T_{2,K}\}$. 
Obviously our angler wants to have as much satisfaction as possible and he has to leave the lake before the fixed moment. Therefore his goal is to find two stopping times $\tau_1^*$ and $\tau_2^*$ such that the expected gain is maximized 
	\begin{equation}
		EZ(\tau_1^*,\tau_2^*)=\sup_{\tau_1 \in \T}\sup_{\tau_2 \in \T^{\tau_1}}EZ(\tau_1,\tau_2),
	\end{equation}
where $\tau_1^*$ corresponds to the moment, when he eventually should change the place and $\tau_2^*$, when he should stop fishing. These stopping times should be less than the fixed moment $t_0$. The process $Z(s,t)$ is piecewise-deterministic and belongs to the class of semi-Markov processes. The optimal stopping of similar processes was studied by ~\cite{bosgou93:semi}.
The special representation of stopping times from $\T$ and $\T^s$ is applied. It allows to use the dynamic programming methods to find these two optimal stopping times and to specify the expected satisfaction of the angler. The way of the solution is similar to the methods used by ~\cite{karsza07:risk}. Let us first observe that by the properties of conditional expectation we have
	\begin{eqnarray*}\label{AK double stop b}
		EZ(\tau_1^*,\tau_2^*)
	 	=\sup_{\tau_1 \in \T}E\{E\left[Z(\tau_1, 				\tau_2^*)|\F_{\tau_1}\right]\} \nonumber 
	 	= \sup_{\tau_1 \in \T}EJ(\tau_1),
	\end{eqnarray*}
where
	\begin{equation}\label{AK one stop}
		J(s)=E\left[Z(s,\tau_2^*)|\F_s\right]=\esssup_{\tau_2 \in \T^{s}}E\left[Z(s,\tau_2)|\F_s\right].
	\end{equation}
Therefore in order to find $\tau_1^*$ and $\tau_2^*$, we have to calculate $J(s)$ first. The process $J(s)$ corresponds to the value of revenue function in the one stopping problem if the observation starts at the moment $s$.
\section{Construction of the optimal second stopping time}
\label{AK second}
In this section, we will find the solution of one stopping problem defined by~(\ref{AK one stop}). We will first solve the problem for fixed number of claims, next we will consider the case with infinite number of claims. In this section we fix $s$ - the moment when the change took place and $m=M_s$ - the mass of the fishes at the time $s$.
\subsection{Fixed number of claims}\label{second fixed}
In this subsection we are looking for optimal stopping time  $\tau_{2,0,K}^*:=\tau_{2,K}^*$ such that
	\begin{equation}
  	E\left[Z(s,\tau_{2,K}^*)|\F_{s}\right]=\esssup_{\tau_{2,K} \in \T^{s}_{0,K}}E\left[Z(s,\tau_{2,K})|\F_{s}\right],
	\end{equation}
where $s\geq 0$ is a fixed time when the position was changed and $K$ is the maximum number of claims which can occur.
Let us define
	\begin{equation}\label{Gamma_{n,K}^s}
  	\Gamma_{n,K}^s=\esssup_{\tau_{2,n,K} \in \T_{n,K}^s}E\left[Z(s,\tau_{2,n,K})|\F_{s,n}\right]
  	=E\left[Z(s,\tau^*_{2,n,K})|\F_{s,n}\right],
	\end{equation}
for $n=K,\dots,1,0$ and observe that $\Gamma_{K,K}^s=Z(s,T_{2,K})$. The crucial role in our subsequent considerations plays the following lemma (see ~\cite{bre81:poiproc}).
	
	\begin{lemma}\label{AK representation}
	If $\tau_1 \in \T$, $\tau_2 \in \T^s$, then there exist \mbox{$R_{1,n}\in \M(\F_n)$} and
	\mbox{$R_{2,n}\in \M(\F_{s,n})$} respectively, such that  
   		$\tau_i \wedge T_{i,n+1}=(T_{i,n}+R_{i,n})\wedge T_{i, n+1}$ on $\{\omega:\tau_i(\omega) \geq T_{i,n}(\omega)\}$, $i \in \{1,2\}$, a.s.
	\end{lemma}  
Now we can derive the dynamic programming equations satisfied by $\Gamma_{n,K}^s$. To simplify notation we write $M_t=M_t^s$ for $t\leq s$, $M_n=M_{T_{1,n}}$, $M_n^s=M_{T_{2,n}}^s$ and $\bar{F_i}=1-F_i$.
	
	\begin{theorem}\label{Gamma recursion}
		Let $s \geq 0$ and $\Gamma_{K,K}^s= Z(s,T_{2,K})$. For $n= K-1, K-2, \dots, 0$ we have
			\begin{eqnarray*}
    		\Gamma_{n,K}^s&=& \esssup_{R_{2,n}\in \M(\F_{s,n})} \vartheta_{n,K}(M_s,s,M_{n}^s,T_{2,n},R_{2,n}) \ a.s.,    		
    	\end{eqnarray*}
    	where	
 			$\vartheta_{n,K}(m,s,\ti m,t,r) = \Ind_{\{t \leq t_0 \}}\bigg\{\bar{F_2}(r)
 																					[\Ind_{\{r \leq t_0-t\}}w_2(m,s,\ti m,t+r)-C\Ind_{\{r>t_0-t\}}]
    																			+ E\left[\Ind_{\{ S_{2,n+1}\leq r \}}\Gamma_{n+1,K}^s|\F_{s,n}\right]\bigg\}
    		-C\Ind_{\{t > t_0 \}}    
  		$.
	\end{theorem}
   It can be shown that there exists $R^*_{2,n}$ such that 
  		$\Gamma_{n,K}^s = \vartheta_{n,K}(M_s,s,M^s_n,T_{2,n},R^*_{2,n})$ for $n \leq K-1$.
	\begin{theorem}\label{tau_2^* representation}
		Let $R_{2,i}^*$ be the sequence of $\F_{s,i}$-measurable random variables $($fix $R_{2,K}^*=0$$)$ and 
		$\eta_{n,K}^s=K\wedge \inf\{i\geq n: R_{2,i}^*<S_{2,i+1}\}$, $n=0,\dots,K$. Then  			
		$\Gamma_{n,K}^s=E\left[Z(s,\tau_{2,n,K}^*)|\F_{s,n}\right]$,
		where $\tau_{2,n,K}^*=T_{2,\eta_{n,K}^s}+R_{2,\eta_{n,K}^{s}}^*$.
	\end{theorem}

	\begin{lemma}\label{AK Gamma 2} 
		$\Gamma_{n,K}^s=\gamma_{K-n}^{s,M_s}(M_{n}^s,T_{2,n})$ for $n=K,\dots,0$, 
		where the sequence of functions	$\gamma_j^{s,m}$ is given recursively as follows:
			\begin{eqnarray}\label{gamma recursive a}
    		\gamma_{0}^{s,m}(\ti m,t)&=&\Ind_{\{t\leq t_0\}}w_2(m,s,\ti m,t)-C\Ind_{\{t > t_0\}},\nonumber \\
    		\gamma_{j}^{s,m}(\ti m,t)&=&\Ind_{\{t\leq t_0\}}\sup_{ r\geq 0}\kappa_{2,\gamma_{j-1}^{s,m}}(m,s,\ti m,t,r)
    																	-  C\Ind_{\{t > t_0\}},
    	\end{eqnarray}
    where		
    	$	\kappa_{2,\delta}(m,s,\ti m,t,r)=
    		\bar{F_2}(r)[\Ind_{\{r \leq t_0-t \}}w_2(m,s,\ti m,t+r)-C\Ind_{\{r > t_0-t\}}]
    		+ \int_0^r dF_2(z) \int_0^\infty \delta(\ti m+x,t+z)dH_2(x)$.
	\end{lemma}
Let us denote $\alpha_i=f_i/\bar{F}_i$ and let us set $\Delta_i(a)=E\left[g_i(a+X_i)-g_i(a)\right]$.
		The sequence of functions $\gamma_j^{s,m}$ can be expressed as follows 
\small			\begin{eqnarray}\label{gamma recursive 2}
    		\gamma_{j}^{s,m}(\ti m,t)&=&\Ind_{\{t\leq t_0\}}\bigg\{ w_2(m,s,\ti m,t)+y_{2,j}(\ti m-m,t-s,t_0-t) \bigg\}
    	-C\Ind_{\{t > t_0\}}
\normalsize		\end{eqnarray} 
		where  $y_{2,j}(a,b,c)$ is given recursively as follows $y_{2,0}(a,b,c)= 0$ and 
		$	y_{2,j}(a,b,c)=\max_{ 0\leq r \leq c}\phi_{2,y_{2,j-1}}(a,b,c,r)$,
 		$$ \phi_{2,\delta}(a,b,c,r)=  \int_0^r \bar F_2(z)
    													\{\alpha_2(z)\left[\Delta_2(a)+E\delta(a+X_2,b+z,c-z)\right]
    							- c_2'(b+z)\}dz.$$
The second optimal stopping time is constructed similarly like in ~\cite{fersie97:risk}. Let $B=B([0,\infty)\times[0,t_0]\times[0,t_0])$ be the space of all bounded, continuous functions with the norm \mbox{$\left\|\delta\right\|=\sup_{a,b,c}|\delta(a,b,c)|$}. It is complete space. Let us define the operator $\Phi_2: B \rightarrow B$ as  
	\begin{equation}\label{AK Phi_2}
		(\Phi_2\delta)(a,b,c)=\max_{0\leq r\leq c}\phi_{2,\delta}(a,b,c,r).
	\end{equation}
We have $y_{2,j}(a,b,c)=(\Phi_2y_{2,j-1})(a,b,c)$ and by (\ref{gamma recursive 2}) there exists a function $r_{2,j}^*(a,b,c)$ such that $y_{2,j}(a,b,c)=\phi_{2,y_{2,j-1}}(a,b,c,r_{2,j}^*(a,b,c))$ and  
		$	\gamma_{j}^{s,m}(\ti m,t)=\Ind_{\{t\leq t_0\}}\bigg\{ w_2(m,s,\ti m,t)+\phi_{2,y_{2,j-1}}(\ti m-m,t-s,t_0-t,r_{2,j}^*(\ti m-m,t-s,t_0-t))\bigg\}
    	-C\Ind_{\{t > t_0\}}$.	
The consequence of the foregoing considerations is the optimal stopping
times $\tau_{2,n,K}^*$ in following form:
	\begin{theorem}\label{AK solution 2}
	Let $R_{2,i}^*  = r_{2,K-i}^*(M_{i}^s-M_s,T_{2,i}-s,t_0-T_{2,i})$ for $i=0,1,\dots,K$ moreover
			$\eta^s_{n,K} = K\wedge \inf \{i\geq n: R_{2,i}^*< S_{2,i+1}\}$, then 
	the stopping time $\tau_{2,n,K}^*=T_{2,\eta^s_{n,K}}+ R_{2,\eta^s_{n,K}}^*$ 
is optimal in the class $\T^s_{n,K}$ and $\Gamma_{n,K}^s=E\left[Z(s,\tau^*_{2,n,K})|\F_{s,n}\right]$.
\end{theorem}
\subsection{Infinite number of claims}
The solution of one stopping problem, related to construction of the second stopping moment, for infinite number of claims is obtained under assumption that $F_2(t_0)<1$.  
	If $F_2(t_0)<1$ then the operator $\Phi_2:B \rightarrow B$ defined by~(\ref{AK Phi_2}) is a contraction.
	By (\ref{gamma recursive 2}) and the contraction properties we get by fixed point theorem that
		there exists $y_2 \in B$ such that $y_2=\Phi_2 y_2$  and $\lim_{K\rightarrow \infty }\|y_{2,K}-y_2\|=0$. 
It implies that $y_2$ is measurable and $\gamma^{s,m}=\lim_{K\rightarrow \infty}\gamma^{s,m}_K$ is given by
	\begin{equation*}\label{gamma^{s,m}}	
			\gamma^{s,m}(\ti m,t )
			=\Ind_{\{ t \leq t_0\}}\left[w_2(m,s,\ti m,t)+y_2(\ti m-m,t-s,t_0-t)\right]-C\Ind_{\{t > t_0\}}.
	\end{equation*}		
We can now calculate the optimal strategy and the expected gain after changing place.
	\begin{theorem}\label{solution 2}
	If $F_2(t_0)<1$ and has the density function $f_2$, then\vspace{-2ex}
		\begin{itemize}
  		\item[\rm(i)] for $n \in \Nat$ the limit $\tau_{2,n}^*=\lim_{K\rightarrow \infty}\tau_{2,n,K}^*\ a.s.$ exists and 
  		$\tau_{2,n}^* \leq t_0$ is an optimal stopping rule in the set $\T^s\cap \{\tau \geq T_{2,n}\}$,
   		\item[\rm(ii)] $ E\left[Z(s, \tau_{2,n}^*)|\F_{s,n}\right]=\gamma^{s,m}(M_{n}^s,T_{2,n})$ a.s.
		\end{itemize}
	\end{theorem}
It can be proved that the function $\gamma^{s,m}(m,s)$ with respect to $s$ is left-hand differentiable. It allows to construct the optimal the  first optimal stopping moment similarly as the second one. To this end we have to take as the payoff function   
			$\gamma^{s,m}(m,s)=\Ind_{\{s\leq t_0\}}u(m,s)-C\Ind_{\{s > t_0\}}$, 
		where $u(m,s)=g_1(m)-c_1(s)+g_2(0)-c_2(0)+\bar y_2(t_0-s)$ is continuous, bounded, measurable with bounded left-hand derivatives with respect to 
		$s$.
The conditional value function of the second optimal stopping problem has the form:
	\begin{equation}\label{J(s)}
		J(s)=E\left[Z(s,\tau_2^*)|\F_s\right]=\gamma^{s,M_s}(M_s,s)\ \ \ a.s.
	\end{equation}	
The first optimal stopping moment in the considered problem is equal $\tau_1^*$ such that $J(\tau_1^*)=\sup_{\tau\in\T}EJ(\tau)$.
\begin{example}
	If $S_{2}$ has exponential distribution with constant hazard rate $\alpha_2$, $g_2$ is increasing and
	concave, $c_2$ is convex and $t_{2,n}=T_{2,n}$, $m_n^s=M_n^s$ then
	\begin{equation*}\label{tau^*_{2,n}}
		\tau^*_{2,n}=\inf \{ t \in \left[t_{2,n},t_0 \right] : \alpha_2[Eg_2(m_n^s + X_2-m)-g_2(m_n^s -m)]\leq c_2'(t-s) \},
	\end{equation*}
	where $s$ is the moment of changing place.
	Moreover if $S_{1}$ has exponential distribution with constant hazard rate $\alpha_1$, $g_1$ is increasing and concave, $c_1$ is convex and $t_{1,n}=T_{1,n}$, $m_n=M_n$ then
	$$\tau^*_{1,n}=\inf \{ s \in \left[t_{1,n},t_0 \right] : \alpha_1\left[Eg_1(m_n+ X_1)-g_1(m_n)\right]\leq  c_1'(s) \}$$
\end{example}
	If for $i=1$ and $i=2$ the functions $g_i$ are increasing and convex, $c_i$ are concave and $S_{i}$ have exponential distribution with constant hazard rate $\alpha_i$ then $\tau^*_{1,n}=\tau^*_{2,n}=t_0$ for $n \in \Nat$.
\section{Conclusions}	
This article presents the solution of double stopping problem in "fishing model" for finite horizon. The analytical properties of the reward function in one stopping problem played the crucial rule in our considerations and allowed us to extend the problem to double stopping.  It is easy to generalize our model and the solution to multiple stopping problem.
\nopagebreak  \vspace{1ex}

\small
\begin{minipage}[h]{0.95\textwidth}

\end{minipage}
\end{document}